\documentclass[12pt,draft]{amsart}

\usepackage{enumerate}

\usepackage{color}

\newtheorem{teo}{Theorem}
\newtheorem{lem}[teo]{Lemma}
\newtheorem{prop}[teo]{Proposition}

\theoremstyle{definition}
\newtheorem{dfn}[teo]{Definition}
\newtheorem{rk}[teo]{Remark}
\newtheorem{ex}[teo]{Example}

\def\<{\langle}
\def\>{\rangle}

\def\b{\beta}
\def\a{\alpha}

\def\e{\varepsilon}

\def\r{\rho}
\def\t{\tau}

\def\M{{\mathcal M}}

\def\t1p{$\operatorname{(T_1P)}$}
\def\T{{\mathbb T}}

\def\Prim{\operatorname{Prim}}

\def\Im{\operatorname{Im}}

\def\ss{\subset}

\newcommand{\wh}[1]{\widehat{#1}}

\begin{document}
\title{A $C^*$-analogue of Kazhdan's property (T)}
\author{
A.~A. Pavlov }
\address{Department of Geography,
Moscow State University, 119 992 Moscow,  Russia}
\email{axpavlov@mail.ru}
\thanks{
Partially supported by the RFBR (grant 05-01-00982) and the Grant
"Universities of Russia"}
\author{
E.~V.~Troitsky }
\address{Dept. of Mech. and Math.,
Moscow State University, 119 992 Moscow,  Russia}
\email{troitsky@mech.math.msu.su} \urladdr{
http://mech.math.msu.su/\~{}troitsky}



\begin{abstract}
This paper deals with a "naive" way of generalization of the
Kazhdan's property (T) to $C^*$-algebras. This approach differs
from the approach of Connes and Jones, which has already
demonstrated its utility. Nevertheless it turned out that our
approach is applicable to the following rather subtle question in
the theory of $C^*$-Hilbert modules.

We prove that a separable unital $C^*$-algebra $A$ has property MI
(\emph{module-infinite}, i.~e. any $C^*$-Hilbert module over $A$
is self-dual if and only if it is finitely generated projective)
if and only if it has not property \t1p (\emph{property {\rm (T)} at one
point}, i.~e. there exists in $\wh A$ a finite dimensional
isolated point).

The commutative case was studied in a previous paper.
\end{abstract}

\maketitle

%%%%%%%%%%%%%%%%%%%%%%%%%%%%%%%%%%%%%%%%%%%%%%%%%%%%%%%%%%%%
 \section{Introduction}
 %\markright{}
%%%%%%%%%%%%%%%%%%%%%%%%%%%%%%%%%%%%%%%%%%%%%%%%%%%%%%%%%%%%
There are at least two basic approaches to the definition of
Kazhdan's property (T) \cite{kazhdan} for topological groups (see
\cite{HarpeValette,BHV,LubotZuk}). The first one uses the notion
of $\varepsilon$-invariant vectors, and property (T) is formulated
in terms of existence of a so-called Kazhdan pair for a
topological group. More precisely, let $G$ be a topological group,
$Q\subset G$ be a compact set, and $\varepsilon>0$. Then the group
$G$ has $T$-property if any unitary representation $(\pi,H)$ of
$G$ which has a $(Q,\varepsilon)$-invariant vector also has a
non-zero invariant vector (in this case the pair $(Q,\varepsilon)$
is called a Kazhdan pair). The second approach is a reformulation
of property (T) in terms of Fell's topology. From this point of
view a topological group $G$ has (T)-property if its identity
representation $1_G$ is isolated in $\mathcal{R}\cup 1_G$, for
every set $\mathcal{R}$ of equivalence classes of unitary
representations of $G$ without non-zero invariant vectors.

An analogue of property (T) for $W^*$-algebras was introduced by
Connes and Jones
 in~\cite{Con1,CJ}. The key notion of this definition was
 the notion of a correspondence which plays a role of a
 representation of a group. Namely, let $A$ and $B$ be von Neumann
 algebras. A correspondence from $A$ to $B$ is a Hilbert space $H$
 which is a left $A$-module and a right $B$-module with commuting
 normal actions. It is possible to introduce some topology on the
 set of all correspondences from $A$ to $B$ which is similar to Fell's
 topology. Let $\mathrm{Id}_A$ be the identity correspondence constructed
in~\cite{Haa}. Then $A$ has (T)-property if there is a
neighborhood $U$ of $\mathrm{Id}_A$ such that any correspondence
in $U$ contains $\mathrm{Id}_A$ as a direct summand. This
definition may be reformulated in terms of central and almost
central vectors. In~\cite{CJ} was proved that a countable discrete
group $G$ with factorial von Neumann algebra $L(G)$ has Kazhdan's
property if and only if $L(G)$ has property (T). This is a very
natural approach and a number of developments and applications was
obtained (in particular, was constructed an example of a
homomorphism $\theta$ of a discrete group $Q$ with property (T)
into ${\rm Out}(\lambda(F\sb \infty)''$) with trivial obstruction
${\rm Ob}\,\theta\in H\sp 3(Q,\T)$ but which cannot be lifted to a
homomorphism from $Q$ to ${\rm Aut}(\lambda(F\sb \infty)''$)
\cite{CJ}).

 For some $C^*$-algebras an analogue of property (T)
 was defined in~\cite{Bekka} by an adaptation of Connes' definition.
 On the whole the approach of B.~Bekka is close
 to the  first approach for topological groups. Namely, let $A$ be
 either a unital $C^*$-algebra admitting a tracial state
  or a finite von Neumann algebra. Then $A$ has Property (T) if
  there exist a finite subset $F$ of $A$ and $\varepsilon>0$ such
  that the following property holds: if a Hilbert $A$-bimodule $H$
  contains a unit $(F,\varepsilon)$-central vector, then $H$ has a
  non-zero central vector. In~\cite{Bekka} it is proved that if
$G$ is a countable discrete group and $A$ a $C^*$-algebra being a
quotient of $C^*_{max}(G)$ such that $C^*_r(G)$ is a quotient of
$A$, then $G$ has property (T) if and only if $A$ has Property (T)
if and only if $L(G)$ has property (T).

In the present paper we discuss the following "naive"
generalization of property (T) to the context of $C^*$-algebras.

\begin{dfn}
A unital $C^*$-algebra $A$ has property \t1p (\emph{property} (T)
(at least) \emph{at one point} of its spectrum) if there exists an
isolated point with respect to the Fell's topology in unitary dual
$\wh A$ such that the corresponding representation is finite
dimensional.
\end{dfn}

\begin{rk}
As it is known (see e.g. \cite[Prop.~14]{HarpeValette},
\cite[Theorem 1.2.6]{BHV}), for locally compact groups the trivial
representation is isolated if and only if all finite dimensional
representations are isolated, if and only if any finite
dimensional representation is isolated. In particular, for
$A=C^*(G)$ our definition is exactly the Kazhdan's one, if $G$ is
a locally compact group.
\end{rk}

We prove in the present paper that property \t1p is responsible
for one very fine property of $C^*$-Hilbert modules over the
algebra under consideration.

\begin{dfn}[\cite{TroPAMS}]\label{dfn:MI}
A unital $C^*$-algebra $A$ is called MI (\emph{module-infinite})
if each countably generated Hilbert $A$-module is projective
finitely generated if and only if it is self-dual. Let us remark
that a projective finitely generated module over a unital algebra
is always self-dual.
\end{dfn}

\medskip\noindent
{\bf Main Theorem.}\ \emph{\label{teo:MInoncommps} Suppose, $A$ is
a separable unital $C^*$-algebra. Then $A$ is {\rm MI} if and only
if $\widehat{A}$ has no isolated point being a finite-dimensional
representation, i.~e. $A$ is not \t1p.}

\medskip
For commutative algebras this was proved in~\cite{TroPAMS}. The
research was motivated by the study of conditional expectations of
finite index related to group  actions
\cite{FMTZAA,TroPAMS,seregin} and \cite[Sect.~4.5]{MaTroBook}.

\begin{rk}\label{rk:equivspecprimfindim}
The property of existence of isolated points can be formulated for
$\wh A$ and for $\Prim(A)$. Fortunately, for finite-dimensional
representations it does not depend on the space. Indeed, a
difference could arise in the case if an isolated point under
consideration has several pre-images. But finite-dimensional
irreducible representations are equivalent if and only if they
have the same kernel  (see e.g. Example 5.6.2 and Theorem 5.6.3 in
\cite{Murphy}).  This will allow us to choice an appropriate
setting to simplify argument.
\end{rk}

\medskip
\noindent {\bf Acknowledgment.\/} The present research is a part of
joint research programm of the second author in Max-Planck-Institut
f\"ur Mathematik (MPI) in Bonn. He would like to thank the
MPI for its kind support and
hospitality while this work has been completed.

%%%%%%%%%%%%%%%%%%%%%%%%%%%%%%%%%%%%%%%%%%%%%%%%%%%%%%%%%%%%
 \section{Preliminaries and reminding}\label{sec:prelim}
 %\markright{}
%%%%%%%%%%%%%%%%%%%%%%%%%%%%%%%%%%%%%%%%%%%%%%%%%%%%%%%%%%%%

The necessary information about Hilbert C*-modules can be found
  in \cite{Lance} and \cite{MaTroBook,MaTroJMS}.
  Let us remind only some facts and notations
  for convenience of the reader.

The standard Hilbert module over a $C^*$-algebra $A$ is denoted by
$l_2(A)$ or $H_A$.

  \begin{prop}[\cite{MaTroBook}, Proposition 2.5.5]\label{prop:funct_on_H_A}
   Consider the set of all sequences $f=(f_i), f_i\in A,
  i\in \mathbb{N}$, such that the partial sums of the series
  $\sum f_i^*f_i$  are uniformly bounded. If $A$ is a unital
  $C^*$-algebra, then this set coincides with $H_A'$, the action of
  $f$ on $H_A$ is defined by the formula
  \begin{gather*}
    f(x)=\sum_{i=1}^\infty f_i^*x_i,
\end{gather*}
where $x=(x_i)\in H_A$, and the norm of the functional $f$ is
satisfies
\begin{gather*}
    \|f\|^2=\sup\limits_N\left\|\sum_{i=1}^\infty f_i^*f_i\right\|.
\end{gather*}
  \end{prop}

\begin{prop}[\cite{MaTroBook}, Theorem 5.1.6]\label{prop:crit_selfd_H_A}
Let $A$ be a $C^*$-algebra. Then the following conditions are
equivalent:

\begin{enumerate}[\rm (i)]
    \item the Hilbert $C^*$-module $H_A$ is self-dual;
    \item the $C^*$-algebra $A$ is finite-dimensional.
\end{enumerate}
\end{prop}

Let us remind also some notions and notations about spectrum of a
$C^*$-algebra $A$ (for more information see~\cite{Dix,Ped}). Let
$\widehat{A}$ be the space of all equivalence classes of
irreducible representations of $A$ and $\mathrm{Prim}(A)$ be the
space of primitive ideals of $A$. For $a\in A$ and $I\in
\mathrm{Prim}(A)$ we will denote by $\|a+I\|$ the norm of element
$a+I$ in the factor algebra $A/I$.

\begin{teo}
[Dauns-Hofmann, \cite{DaunsHofmann,EllOle}]
\label{teo:Dauns-Hofmann} Let $A$ be a $C^*$-algebra, let $x$ be
an element of $A$, and let $f$ be a bounded continuous
scalar-valued function on $\mathrm{Prim}(A)$, the space of
primitive ideals of $A$ endowed with the hull-kernel topology.
Then there exists a unique element $fx$ of $A$ such that
\begin{gather*}
    (\widehat{fx})(I)=f(I)\widehat{x}(I),\qquad
    I\in\mathrm{Prim}(A),
\end{gather*}
where $\widehat{x}(I)=x+I\in A/I$.
\end{teo}

A point $S\in\mathrm{Prim}(A)$ is called a \emph{separated point}
if for any $I\in\mathrm{Prim}(A)$ that is not an accumulation
point of $S$ there are disjoint neighborhoods of $S$ and $I$.

\begin{teo}[\cite{DixASM61}]
Suppose $A$ is a $C^*$-algebra. Then

\begin{enumerate}[\rm (i)]
    \item a point $I\in \mathrm{Prim}(A)$ is separated if and
only if for any $x\in A$ the function $I\mapsto \|x+I\|$ is
continuous at $I$;
    \item if $A$ is separable, then the set of separated
points of $\mathrm{Prim}(A)$ is $G_\delta$ and it is dense in
$\mathrm{Prim}(A)$.
\end{enumerate}
\end{teo}

%%%%%%%%%%%%%%%%%%%%%%%%%%%%%%%%%%%%%%%%%%%%%%%%%%%%%%%%
\section{Some properties of $\mathrm{Prim}(A)$.}
%%%%%%%%%%%%%%%%%%%%%%%%%%%%%%%%%%%%%%%%%%%%%%%%%%%%%%%%
Suppose, $A$ is a unital $C^*$-algebra.

\begin{lem}\label{lem:Urysohn}
Let $U$ be a neighborhood of a point $\rho\in \widehat{A}$. Then
there exists a positive element $a\in A$ of norm $1$ such that
\begin{enumerate}[\rm(i)]
    \item $\|\rho (a)\|=1,$
    \item $\pi (a)=0$ for all $\pi\in\widehat{A}\setminus U$.
\end{enumerate}
The same is true in $\Prim(A)$.
\end{lem}

\begin{proof}
By~\cite[Lemma 3.3.3]{Dix} there exists a positive element $x\in
A$ such that the set $Z=\{\pi\in\widehat{A} : \|\pi(x)\|>1\}$
contains $\rho$ and belongs to $U$.  Let $\|\rho(x)\|=t$, $1<t\le
\|x\|$, and let $t_1\in (1,t)$. Let $f$ be a continuous function
on $[0,\|x\|]$, equal to 0 on $[0,t_1]$, equal to 1 on $[t,\|x\|]$
and linear on $[t_1,t]$, and let $a=f(x)$. Then $a\ge 0$,
$\|a\|=1$. Besides,
\begin{gather*}
    \|\rho(a)\|=\|\rho(f(x))\|=\|f(\rho(x))\|=1
\end{gather*}
and
\begin{gather*}
    \|\pi(a)\|=\|\pi(f(x))\|=\|f(\pi(x))\|=0
\end{gather*}
for all $\pi\in\widehat{A}\setminus U$, because $\|\pi(x)\|\le
1<t_1$.
\end{proof}

\begin{lem}\label{lem:sep_points2}
Let $\{S_i\}$ be a sequence of different separated points of
$\mathrm{Prim}(A)$. Then there exist a subsequence $\{S_{i(j)}\}$
of $\{S_i\}$ and  neighborhoods $V_{j}$ of $S_{i(j)}$, which do not contain
the remaining points of the subsequence.
\end{lem}

\begin{proof}
Suppose, $\{S_i\}$ has no accumulation points among its members.
Then each of them has a neighborhood with no other points of the
sequence.

Now, suppose that $S_k$ is an accumulation point. Without loss of
generality (we can pass to a subsequence) we assume that $S_k$ is
a limit point. Let $V'_1$ and $V_1$ be disjoint neighborhoods of
$S_k$ and $S_1=:S_{i(1)}$ (since $\Prim(A)$ is a $T_0$ space,
separated points there form a Hausdorff space, in particular, any
finite number of them has pairwise disjoint neighborhoods). Then
there exist $S_{i(2)}\in V'_1$, $S_{i(2)}\ne S_k$, and their
disjoint neighborhoods $W'_2 \ni S_k$ and $W_2\ni S_{i(2)}$. Take
$V'_2:=W'_2\cap V'_1$ and $V_2:=W_2\cap V'_1$. And so on.
\end{proof}

%%%%%%%%%%%%%%%%%%%%%%%%%%%%%%%%%%%%%%%%%%%%%%%%%%%%%%%
\section{The property DINC}
%%%%%%%%%%%%%%%%%%%%%%%%%%%%%%%%%%%%%%%%%%%%%%%%%%%%%%%
\begin{dfn}\label{dfn:DIMNCps}
A unital C*-algebra is said to be DINC (\emph{divisible infinite
for the noncommutative case}) if for any sequence $u_i\in A$ of
elements of norm $1\ge \|u_i\|\ge C>0$   there exist a subsequence
$i(k)$ and elements $b_k \in A$ of norm $1$ such that

(i) the partial sums of the series $\sum_k b_k^* b_k$ are
uniformly bounded, and

(ii) for each $k$
\begin{equation}\label{eq:proprhok}
  \|u_{i(k)}b_k\|\ge C/2.
\end{equation}
\end{dfn}

The following result is a generalization of~\cite[Theorem
32]{TroPAMS}.

\begin{teo}\label{DIMtoMIncps}
If a unital $C^*$-algebra $A$ is {\rm DINC}, then it is {\rm MI}.
\end{teo}

\begin{proof}
We have to prove that if a countable generated Hilbert $A$-module
$\M$ is not finitely generated projective, then it is not
self-dual. By the Kasparov stabilization theorem \cite{KaspJO}
(see also \cite{Lance,MaTroBook}) one has $\M\oplus l_2(A)\cong
l_2(A)$.
 Denote by $p_\M:l_2(A)\to\M\subset l_2(A)$ the
corresponding orthoprojection. Let $p_j: l_2(A)\to E_j \cong A^j$
be the orthoprojection on the first $j$ standard summands of
$l_2(A)$ and $q_j$ the orthoprojection on the $j$-th standard
summand in such a way that $p_j=q_1+\dots +q_j$. Two possibilities
can arise: 1) $\|(1-p_j) p_\M\|\to 0$ as $j\to\infty$, and 2) the
opposite case.

1) Let us show that in this case $\M$ is finitely generated
projective, i.e. this case is impossible under our assumptions.
One can argue as in \cite{MF}: for a sufficiently large $j$ the
operator
 $$
J(x)=\left\{
 \begin{array}{ll}
p_j(x) & \mbox{ if } x\in \M,\\
x & \mbox{ if } x\in \M^\bot\cong l_2(A)
 \end{array}
\right.
 $$
is close to identity, hence an isomorphism. It maps $\M$
isomorphically onto a direct summand of $E_j$.

2) In this case consider the matrix of the orthogonal projection
$p_\M$. This is an adjointable (in fact self-adjoint) operator
$p_\M:l_2(A)\to l_2(A)$. Hence $\|p_j p_\M (1-p_k)\|\to 0$ as
$k\to\infty$ for fixed $j$. Indeed, it is sufficient to verify
that $\|q_m p_\M (1-p_k)\|\to 0$ for each $m=1,\dots, j$, i.e.
$$
\sup_{\|x\|=1} |\<e_m, p_\M (1-p_k) x\>|= \sup_{\|x\|=1}
|\<(1-p_k) p_\M e_m, x\>|= \| (1-p_k) p_\M e_m\|\to 0.
$$
The last is evident because $p_\M e_m$ is a fixed element of
$l_2(A)$. So, we have in this case the following relations:
$$
\|(1-p_j) p_\M\|\not\to 0\quad (j\to\infty), \qquad \|p_j p_\M
(1-p_k)\|\to 0\quad (k\to\infty,\: j\mbox{ is fixed}).
$$
Using these two observations one can choose decompositions
$l_2(A)=M_1\oplus M_2 \oplus \dots$ of the domain and
$l_2(A)=N_1\oplus N_2 \oplus \dots$ of the range of $p_\M$ in such
a way that
\begin{itemize}
\item $M_i$ and $N_i$ are free modules generated by several consequent
vectors
$$e_{\mu(1,i)},\dots, e_{\mu(m_i,i)} \qquad \mbox{and}
\qquad e_{\nu(1,i)},\dots, e_{\nu(n_i,i)}
$$
of the standard base;
\item there exist elements $v_i\in M_i$, $\|v_i\|\le 1$, such that
the projection of $p_\M(v_i)$ onto $N_i$ is not small, more
precisely
$$
\|(p_{\nu(n_i,i)}-p_{\nu(1,i)-1})p_\M(v_i)\|\ge C
$$
for some fixed $C>0$ for any $i$;
\item for these elements the projection of $p_\M(v_i)$ onto
$N_1\oplus \dots \oplus N_{i-1} \oplus N_{i+1}\oplus \dots $ is
small, more precisely
$$
\|(1-p_{\nu(n_i,i)}+p_{\nu(1,i)-1})p_\M(v_i)\| < \frac \e{2^i}
$$
for any beforehand fixed sufficiently small $\e>0$ and any $i$.
\end{itemize}
We find these elements by induction over $i$. By the supposition
there exists a number $K>0$ such that for any $j$
\begin{equation}\label{eq:nestrem}
    \|(1-p_j) p_\M w_j\| \ge K
\end{equation}
for some $w_j$ of norm 1. Evidently, $K\le 1$. Let $\e= K/4$. Let
$x$ be any vector of norm 1 from $\M$. Then up to $\e/4$ the
vector $x$ is in $M_1=N_1=E_s$ for some
$s=\mu(m_1,1)=\nu(n_1,1)=m_1=n_1$ and the conditions hold with
$v_1$ being the projection of $x$ onto $M_1$ and $C=1-\e \ge K -
\e$. Now choose a number $t> \mu(m_1,1)$ such that $\|p_s p_\M
(1-p_t)\|< \e/16$ and after that a number $d>\nu(n_1,1)$ such that
$ \|(1-p_d) p_\M p_t \| < \e/16$. Choose $w_d$ as in
(\ref{eq:nestrem}) and a number $r>t$ such that
$\|(1-p_r)(1-p_t)w_d\|< \e/16$. Then
\begin{multline*}
\|(1-p_d) p_\M p_r (1-p_t) w_d\|\ge \|(1-p_d) p_\M (1-p_t) w_d\|-
\frac \e{16}
\\
\ge \|(1-p_d) p_\M w_d\|- \frac \e{8}\ge K - \frac \e{8}.
\end{multline*}
Now choose a number $l>d$ such that
$$
\|p_l(1-p_d) p_\M p_r (1-p_t) w_d\|\ge K -\frac \e{4}
$$
and
\begin{equation}\label{eq:ozenhvosta}
    \|(1-p_l) p_\M p_r (1-p_t) w_d\|< \frac \e{8}.
\end{equation}
Let $\mu(m_2,2):=r$, $\nu(n_2,2):=l$, $v_2:= p_r (1-p_t) w_d$.
Then $\|v_2\|\le 1$ and $v_2\in M_2$, because $t>\mu(m_1,1)$.
Since $l>d>\nu(n_1,1)$, one has
\begin{multline*}
\|(p_{\nu(n_2,2)}-p_{\nu(n_1,1)})p_\M(v_2)\| \ge\|(p_l -p_d)p_r
(1-p_t) w_d\|
\\
=\|p_l(1 -p_d)p_r (1-p_t) w_d\| \ge K-\frac \e{4}.
\end{multline*}
Since $\|p_s p_\M (1-p_t)\|< \e/16$, one has
\begin{multline*}
\e/16 > \|p_{\nu(n_1,1)} p_\M (1-p_t) p_r  w_d\| =
\|p_{\nu(n_1,1)} p_\M  p_r (1-p_t) w_d\| = \|p_{\nu(n_1,1)} p_\M
v_2\|.
\end{multline*}
Together with (\ref{eq:ozenhvosta}) this gives the last necessary
property. Proceeding in such a way we obtain the desired
decompositions and elements with
$$
C=K-\e\ge \frac 34\cdot K.
$$
Denote $u_i:=p_\M(v_i)$. Let $z_i$ be the orthoprojection of $u_i$
onto $N_i$. Then $C\le \|u_i\|\le 1$, $\|z_i-u_i\|<\e/2^i$,
$\|z_i\|\le (1+\e/2^i)$.

According to Definition \ref{dfn:DIMNCps} let us choose elements
$b_k$ for $\<z_k,z_k\>$ (for the sake of notational brevity  we
assume that we do not need to pass to a subsequence the  second
time).  Then the formula
\begin{equation}\label{ur:deffunkbsps}
\b(x)=\sum_k   b_k^*  \<z_k,x\>
\end{equation}
defines an $A$-functional on $l_2(A)$. By
Proposition~\ref{prop:funct_on_H_A} to verify this, it is
sufficient to prove that partial sums of the series
$$
\sum_i \b(e_i)\b(e_i)^*
$$
are uniformly bounded. If $e_i \in N_m$, then
$$
\b(e_i)=\sum_{k\ne m} b_k^* \<z_k,e_i\> + b_m^* \<z_m, e_i\>=
b_m^* \<z_m, e_i\>.
$$
Hence,
\begin{eqnarray*}
  \sum_i \b(e_i)\b(e_i)^* &=& \sum_m \sum_{e_i\in N_m}\left(b_m^*
  \<z_m, e_i\>\right)
\left(b_m^* \<z_m, e_i\>\right)^* \\
    &=& \sum_m b_m^* \left(\sum_{e_i\in N_m}\<z_m, e_i\>\<z_m, e_i\>^*
    \right) b_m
    \\
    &=& \sum_m b_m^* \<z_m,z_m\> b_m \le (1+\e)^2 \sum_m b_m^* b_m.
\end{eqnarray*}

 Let us show
that there is no adjointable functional $\a$ on $l_2(A)$ such that
 $\a|_\M=\b|_\M$, and hence, $\b|_\M$ is a non-adjointable functional on $\M$
and $\M$ is not a self-dual module. Indeed, suppose, there exists
an element $a=(a_1,a_2,\dots)\in \M\ss l_2(A)$ such that
$\a(x):=\sum_i a_i x^i= \b(x)$ for any $x\in\M$. Then
$$
\a(z_j)\to 0,\quad [\a(u_j)-\a(z_j)]\to 0,\quad \a(u_j)=\b(u_j),
$$
$$
 [\b(u_j)-\b(z_j)]\to 0, \quad
\b(z_j)=b_j^*\<z_j,z_j\> \not\to 0.
$$
We obtain a contradiction.
\end{proof}

%%%%%%%%%%%%%%%%%%%%%%%%%%%%%%%%%%%%%%%%%%%%%%%%%%%%%%%%%%%%
 \section{A description of spectra of $MI$-algebras}
 %\markright{}
%%%%%%%%%%%%%%%%%%%%%%%%%%%%%%%%%%%%%%%%%%%%%%%%%%%%%%%%%%%%

The following theorem is proved in \cite{FrankZAA} in a much more
generality. We present here a more elementary argument to make the
present text more self-contained.

\begin{teo}\label{teo:finitedimequivconv}
Let a $C^*$-algebra $A$ be an irreducible subalgebra of $B(H)$.
Then the following properties are equivalent:
\begin{itemize}
\item $H$ is infinite-dimensional;
\item there exists a sequence $\a_i \in A$ such that
\begin{enumerate}[\rm (i)]
    \item $\|\a_i\|=1$;
    \item $0\le \a_i \le 1$;
    \item $\sum_{i=1}^n \a_i \le 1$ for any $n$.
\end{enumerate}
\end{itemize}
\end{teo}

\begin{proof}
If $H$ is finite-dimensional, evidently a sequence with the
mentioned properties does not exist.

Now let $H$ be infinite-dimensional. If there exists an element
$\a\ge 0$ of $A$ with infinite spectrum, then we can construct the
desired $\a_i$ with the help of functional calculus on the
spectrum of $\a$, i.e. to find them inside the commutative
$C^*$-algebra generated by $\a$.

So, suppose, that any positive element of $A$ has a finite
spectrum. In this situation any projection $p\in A$ can be
decomposed into a sum of 2 non-trivial projections $p=p_0+p_1$
supposing that the image of $p$ has dimension at least 2. Indeed,
let $h_0$ and $h_1$ be two unit vectors from the image of $p$,
which are orthogonal to each other. By Kadison transitivity
theorem there exists a self-adjoint $\a\in A$ such that
$\a(h_0)=0$ and $\a(h_1)=h_1$. Then $p\a^*\a p=p\a^2 p$ has 0 and
1 as eigenvalues with eigenvectors in $\Im p$. Hence, its spectral
decomposition over $\Im p$ is non-trivial, while the projections
are in $A$ since the spectrum is finite.

Now we can repeat this taking of decomposition $p=p_0+p_1$
infinitely many times starting from some spectral projection and
obtain a sequence of mutually orthogonal projections, which can be
taken as $\a_i$.
\end{proof}

We need the following lemma for the proof of the next theorem.

\begin{lem}\label{lem:17ps}
Suppose, $\widehat{A}$ has an isolated point being a
finite-dimensional representation. Then $A$ is not {\rm MI}.
\end{lem}

\begin{proof}
Let $M\in \mathrm{Prim}(A)$ be the kernel of this isolated
representation $\rho$. Then $M$ is a maximal ideal.  Let
$$ \chi(I)=\left\{
\begin{aligned}
1,\qquad I=M,\\
0,\qquad I\neq M
\end{aligned}\right.
$$
be the characteristic function of $\{M\}$. The set $\{M\}$ is open
and closed, therefore the function $\chi$ is continuous. Let us
consider the finite-dimensional matrix $C^*$-algebra
$M_n=\rho(A)\cong A/M$. By the Dauns-Hofmann theorem
(Theorem~\ref{teo:Dauns-Hofmann}) for any $a\in A$ there exists a
unique element $(\chi a)\in A$ such that $\chi
\widehat{a}=\widehat{\chi a}$. Here\begin{gather*}
    \widehat{a}(I)=a+I, \quad I\in \mathrm{Prim}(A).
\end{gather*}
Suppose,
\begin{gather*}
A_M=\{a\in A : \widehat{a}(I)=0\quad\text{for all}\quad I\neq M\}
\end{gather*}
and $\pi$ is the composition
\begin{gather*}
    A\longrightarrow A/M\stackrel{\cong}{\longrightarrow} M_n.
\end{gather*}
Then by the Dauns-Hofmann theorem the restriction of $\pi$ to
$A_M$ is an isomorphism. Also, $A_M$ is a direct summand in $A$,
because the map $a\mapsto \chi a$ is a projection. Thus $A$ is not
MI (see~Proposition~\ref{prop:crit_selfd_H_A}).
\end{proof}

%\begin{teo}\label{teo:MInoncommps}
%Suppose, $A$ is a separable unital $C^*$-algebra. Then $A$ is {\rm MI} if
%and only if $\widehat{A}$ has no isolated point being a
%finite-dimensional representation, i.~e. $A$ is not \t1p.
%\end{teo}

\begin{proof}[Proof of the Main Theorem]
In one direction this is just Lemma \ref{lem:17ps}.

Now, suppose, $\widehat{A}$ has no isolated point being a
finite-dimensional representation. We will show that $A$ is MI by
demonstrating that $A$ is DINC. In accordance with Remark
\ref{rk:equivspecprimfindim} we will work with $\Prim(A)$.

 We take any sequence $u_i\in A$ of norm
$1\ge\|u_i\|\ge C$. The functions $f_i(I)=\|u_i+I\|$ are lower
semi-continuous on $\mathrm{Prim}(A)$, therefore the sets
$G_i=f_i^{-1}(2C/3,\infty)$ are open. Let us choose separated
points $S_i$ from $G_i$, so $1\ge\|u_i+S_i\|\ge 2C/3$. Passing to
a subsequence if necessary, we can assume  that either
(\textbf{I}) $S_i=S$ for all $i$ or (\textbf{II}) all $S_i$'s are
different.

In the case (\textbf{I}) we can suppose that $S$ is an isolated
point being a kernel of infinite-dimensional representation $\rho
: A\rightarrow B(H)$. Indeed, in the opposite case obviously there
exists a sequence $T_i$ of distinct separated points such that
$T_i\rightarrow S$. Passing from $S$ to $T_i$, we obtain the case
$\textbf{(II)}$.

So, let us study the case (\textbf{I}) with this
supposition on $S$ to be isolated. The set $\{S\}$ is open,
therefore by Lemma~\ref{lem:Urysohn} there exists a positive
element $a\in A$ of norm 1 such that $\|a+S\|=1$ and $\|a+I\|=0$
for all primitive ideals $I\neq S$. By
Theorem~\ref{teo:finitedimequivconv} there exists a sequence
$a_i\in A$ such that
\begin{enumerate}
    \item $\|a_i+S\|=1$;
    \item $0\le a_i+S \le 1+S$;
    \item $\sum_{i=1}^n a_i+S \le 1+S$ for any $n$.
\end{enumerate}
Consider an irreducible representation $\rho$ of $A$ on $H$ with
kernel $S$. Now we choose $x_i, y\in H$ of norm 1 such that
\begin{gather*}
    \|\r(a_i)x_i\|>\|\r(a_i)\|-\varepsilon=1-\e,\qquad
    \|\r(a)y\|>\|\r(a)\|-\varepsilon>1-\e,
\end{gather*}
where $(1-\e)^2> 1/\sqrt{2}$. By the Kadison transitivity theorem
there exist unitary operators $\r(c_i)\in \r(A)$, $\|c_i\|=1$
(cf.~\cite[Theorem 2.8.3.]{Dix}), such that
\begin{gather*}
    y=\r(c_i)\left(\frac{\r(a_i)x_i}{\|\r(a_i)x_i\|}\right).
\end{gather*}
Hence,
\begin{gather*}
    \|a(c_ia_i)+S\|\ge\|\r(a c_i
    a_i)x_i\|=\|\r(a)y\|\|\r(a_i)x_i\|>
    (\|\r(a_i)\|-\varepsilon)
    (\|\r(a)\|-\varepsilon)>1/\sqrt{2}.
\end{gather*}
Put $v_i=ac_ia_i$ and $b_i=v_i^*v_i$. We claim that
\begin{enumerate}
    \item $\|b_i+S\|>1/2$;
    \item $0\le b_i+S \le 1+S$;
    \item $\sum_{i=1}^n b_i+S \le 1+S$ for any $n$;
    \item $\|b_i+I\|=0$ for all primitive ideals $I\neq S$.
\end{enumerate}
The third property follows from the estimation
\begin{gather*}
    \sum_{i=1}^n b_i+S=\sum_{i=1}^n a_i^*c_i^* a^*ac_ia_i+S\le
    \sum_{i=1}^n a_i^*a_i+S\le 1+S.
\end{gather*}
The others properties above are clear. From
$\|b_i\|=\sup\{\|b_i+I\| : I\in\mathrm{Prim}(A)\}$ it follows that
$\|b_i\|=\|b_i+S\|>1/2$ for all $i$ and $\|\sum_{i=1}^n b_i\|\le
1$ for all $n$. Therefore
$$ u_ib_i+I=\left\{
\begin{aligned}
u_ib_i+S,\qquad I=S,\\
0,\qquad I\neq S,
\end{aligned}\right.
$$
and $\|u_ib_i\|=\|u_ib_i+S\|$. Now we choose $\xi_i, \eta_i\in H$
of norm 1 such that
\begin{gather*}
    \|\r(b_i)\xi_i\|>\|\r(b_i)\|-\varepsilon>1/2-\e,\qquad
    \|\r(u_i)\eta_i\|>\|\r(u_i)\|-\varepsilon\ge C/2-\e,
\end{gather*}
where $(1/2-\e)(C/2-\e)> C/4$. By the Kadison transitivity theorem
there exist unitary operators $\r(c_i)\in \r(A)$, $\|c_i\|=1$,
such that
\begin{gather*}
    \eta_i=\r(c_i)\left(\frac{\r(b_i)\xi_i}{\|\r(b_i)\xi_i\|}\right).
\end{gather*}
Hence,
\begin{gather*}
    \|u_i(c_ib_i)\|\ge\|\r(u_i c_i
    b_i)\xi_i\|=\|\r(u_i)\eta_i\|\|\r(b_i)\xi_i\|>
    (\|\r(b_i)\|-\varepsilon)
    (\|\r(u_i)\|-\varepsilon)>C/4.
\end{gather*}
 Thus $c_ib_i$ can serve as a sequence for $u_i$ denoted by $b_i$ in
 Definition~\ref{dfn:DIMNCps}.
 Therefore $A$ is DINC and
 by Theorem~\ref{DIMtoMIncps} it is MI.

So, let us consider the case $\textbf{(II)}$. By
Lemma~\ref{lem:sep_points2} there are neighborhoods $V_i$ of $S_i$
such that $S_j\notin V_i$ if $j\neq i$ (for the sake of notational
brevity we do not pass to a subsequence). Now by
Lemma~\ref{lem:Urysohn} we can find a sequence $y_n$ of positive
elements of $A$ such that
\begin{equation}\label{eq:propela}
  \|y_n\|=1,\qquad \|y_n + S_n\|=1,\quad
    \|{y}_n+ S_m\|=0\qquad\mbox{for all}\qquad m\neq n.
\end{equation}
Now we will pass from the sequence $y_n$ to another sequence $z_n$
with the following properties:
\begin{description}
  \item[(a)] the partial sums of $\sum_k
z_k^*z_k$ are uniformly bounded by $1$,
 \item[(b)] $\|z_i+S_j\|=0$ for all $j\neq i$,
 \item[(c)] $\|z_i+S_i\|=1$  for any $i$.
\end{description}
For this purpose let us define $z_1:=y_1$ and by induction
\begin{equation}\label{eq:defofbi}
  z_{i+1}:=y_{i+1}\left(1-\sum_{k=1}^{i}z_k^*z_k\right)^{1/2}.
\end{equation}
From the item (a) for $z_i$ it follows that $z_{i+1}$ is
well-defined. The item (b) is evident. The item (a) follows from
the estimation:
$$
\sum_{k=1}^{i+1}z_k^*z_k=\sum_{k=1}^{i}z_k^*z_k+
\left(1-\sum_{k=1}^{i}z_k^*z_k\right)^{1/2} (y_{i+1})^2
\left(1-\sum_{k=1}^{i}z_k^*z_k\right)^{1/2} \le 1.
$$
Besides, by (b)
$$
\left(1-\sum_{k=1}^{i-1}z_k^*z_k\right)^{1/2}+S_i=1.
$$
Hence,
$$
\|z_i+S_i\|= \left\|y_{i}+S_{i}\right\|=1
$$
and (c) holds as well.

Now let $S_i=\mathrm{ker}(\rho_i)$ and $H_i$ be the space of the
irreducible representation $\rho_i$. We can choose $\xi_i,
\eta_i\in H_i$ of norm 1 such that 
\begin{gather*}
    \|\r_i(z_i)\xi_i\|>\|\r_i(z_i)\|-\varepsilon=1-\varepsilon,\qquad
    \|\r_i(u_i)\eta_i\|>\|\r_i(u_i)\|-\varepsilon\ge 2C/3-\e,
\end{gather*}
where $(1-\e)(2C/3-\e)\ge C/2$. By the Kadison transitivity
theorem there exist unitary operators $\r_i(c_i)\in \r_i(A)$,
$\|c_i\|=1$, such that
\begin{gather*}
    \eta_i=\r_i(c_i)\left(\frac{\r_i(z_i)\xi_i}{\|\r_i(z_i)\xi_i\|}\right).
\end{gather*}
Thus,
\begin{multline*}
    \|u_i(c_iz_i)\|\ge \|\r_i(u_i(c_iz_i))\|\ge \|\r_i(u_i c_i
    z_i)\xi_i\|
 \\   =\|\r_i(u_i)\eta_i\|\|\r_i(z_i)\xi_i\|>(1-\varepsilon)
    (\|\r_i(u_i)\|-\varepsilon)\ge C/2.
\end{multline*}
 Thus $b_i:=c_iz_i$ is a sequence for $u_i$ as in
 Definition~\ref{dfn:DIMNCps}.
   Therefore $A$ is DINC and
 by Theorem~\ref{DIMtoMIncps} it is MI.
\end{proof}

 Let us consider a couple of examples.

\begin{ex} Let $K(H)$ be the $C^*$-algebra  of compact
operators in an infinite-dimension\-al separable Hilbert space
$H$. For $x, y, z\in H$ suppose $\theta_{x,\,y}(z)=x\< y,z\>$, so
$\theta_{x,\,y}\in~K(H)$. Let $A=\widetilde{K(H)}$ be the algebra
$K(H)$ with an adjoint unit.
 Then $\mathrm{Prim}(A)=\{0, K(H)\}$ and
   all open sets in  $\mathrm{Prim}(A)$ are $\emptyset, \{0\}, \{0,K(H)\}$.
Thus $\{0\}$ is an isolated point corresponding to the identity
(infinite-dimensional irreducible) representation in $H$. Suppose
$\{e_i\}_{i=1}^\infty$ is an orthonormal basis for $H$ and $p_i$
are projections from $H$ onto $\mathrm{span}\{e_i\}$. Then
\begin{enumerate}
    \item $\|p_i\|=1$;
    \item $0\le p_i \le 1$;
    \item $\sum_{i=1}^n p_i \le 1$ for any $n$;
    \item $\|p_i+K(H)\|=0$ for any $i$.
    \end{enumerate}
 Now we take any sequence $u_i\in A$ of norm
$1\ge\|u_i\|\ge C$  and choose $x_i\in H$ of norm 1 such that
\begin{gather*}
    \|u_i x_i\|>\|u_i\|-\varepsilon\ge C-\e>C/2.
\end{gather*}
Put $c_i=\theta_{x_i,\,e_i}$. Then $c_i(e_i)=x_i$ and
\begin{gather*}
    \|u_i(c_ip_i)\|\ge\|u_i c_i
    p_i e_i\|=\|u_i x_i\|> C/2.
\end{gather*}
 Thus $c_ip_i$ can serve as a sequence for $u_i$ denoted by $b_i$ in
 Definition~\ref{dfn:DIMNCps}.
 Therefore $A$ is DINC and
 by Theorem~\ref{DIMtoMIncps} it is MI.
\end{ex}

\begin{ex} Let $\mathbb{A}$ be the Toeplitz algebra and $H^2$ be
the Hardy space. The identity representation of $\mathbb{A}$ in
$H^2$ is  irreducible (cf.~\cite[Theorem 3.5.5]{Murphy}). We claim
that any non zero ideal $I\subset \mathbb{A}$ contains $K(H^2)$.
Indeed, let $u\in I$, $u\neq 0$. Then there exists $x\in H^2$ of
norm 1 such that $u(x)\neq 0$. For any $y\in H^2$ put
\begin{gather*}
p=\theta_{y,\,y},\qquad v=\frac{\theta_{y,\,u(x)}}{\|u(x)\|^2}.
\end{gather*}
Then $vu(x)=y$ and $p=vu\theta_{x,\,x}u^*v^*\in I$. Therefore
$K(H^2)\subset I$.
 Thus the zero ideal
is a unique isolated point in $\mathrm{Prim}(\mathbb{A})$.
Similarly to the previous example one can demonstrate that the
Toeplitz algebra $\mathbb{A}$ is DINC and
 by Theorem~\ref{DIMtoMIncps} it is MI.
\end{ex}

%\bibliographystyle{plain}
%\bibliographystyle{amsplain}
%\bibliography{act03eng}
\def\cprime{$'$}
\providecommand{\bysame}{\leavevmode\hbox to3em{\hrulefill}\thinspace}
\providecommand{\MR}{\relax\ifhmode\unskip\space\fi MR }
\providecommand{\MRhref}[2]{%
  \href{http://www.ams.org/mathscinet-getitem?mr=#1}{#2}
}
\providecommand{\href}[2]{#2}

\end{document}